\documentclass[12pt]{amuc}
\usepackage{amsfonts}
\usepackage{epsfig}
\usepackage{epsfig}
\def\Z{{\mathbb Z}}        
		
\def\F{{\mathcal F}}
\def\H{{\mathcal H}}
\def\E{{\mathbb E}}
\def\P{{\mathbb P}}
\def\L{{\mathcal L}}
\def\1{{\left< 1 \right>}}
\def\0{{\left< 0 \right>}}

\def \MM{{\mathcal M}}

\def \PP{{\mathcal P}}

\def\ascv{\stackrel{\scriptscriptstyle a.s.}{\longrightarrow}}     

\def\mydots{\mathinner{\mkern1mu\raise7pt\vbox{\kern7pt\hbox{.}}\mkern2mu\raise4pt\hbox{.}\mkern2mu\raise1pt\hbox{.}\mkern1mu}}

\def\dude{\mathinner{\mkern1mu\raise8pt\vbox{\kern7pt\hbox{$\scriptscriptstyle\square$}}\mkern2mu\raise4pt\hbox{$\scriptscriptstyle\square$}\mkern2mu\raise0pt\hbox{$\scriptscriptstyle\square$}\mkern1mu}}

\def\RTO{\mathinner{\mkern1mu\raise6pt\vbox{\kern7pt\hbox{$\scriptscriptstyle\square$}}\mkern-2mu\raise2pt\hbox{$\scriptscriptstyle\square$}\mkern-2mu\raise-2pt\hbox{$\scriptscriptstyle\square$}\mkern1mu}}

\def\OTR{\mathinner{\mkern1mu\raise-2pt\vbox{\kern7pt\hbox{$\scriptscriptstyle\square$}}\mkern-2mu\raise2pt\hbox{$\scriptscriptstyle\square$}\mkern-2mu\raise6pt\hbox{$\scriptscriptstyle\square$}\mkern1mu}}

\def\TOT{\mathinner{\mkern1mu\raise4pt\vbox{\kern7pt\hbox{$\scriptscriptstyle\square$}}\mkern-2mu\raise0pt\hbox{$\scriptscriptstyle\square$}\mkern-2mu\raise4pt\hbox{$\scriptscriptstyle\square$}\mkern1mu}}

\def\OTUR{\mathinner{\mkern1mu\raise-2pt\vbox{\kern7pt\hbox{$\scriptscriptstyle\square$}}\mkern-2mu\raise2pt\hbox{$\scriptscriptstyle\square$}\mkern-9.75mu\raise6pt\hbox{$\scriptscriptstyle\square$}\mkern1mu}}

\def\OT{\mathinner{\mkern1mu\raise0pt\vbox{\kern7pt\hbox{$\scriptscriptstyle\square$}}\mkern-2mu\raise4pt\hbox{$\scriptscriptstyle\square$}\mkern1mu}}

\def\TO{\mathinner{\mkern1mu\raise4pt\vbox{\kern7pt\hbox{$\scriptscriptstyle\square$}}\mkern-2mu\raise0pt\hbox{$\scriptscriptstyle\square$}\mkern1mu}}

\def\OTRF{\mathinner{\mkern1mu\raise-4pt\vbox{\kern7pt\hbox{$\scriptscriptstyle\square$}}\mkern-2mu\raise0pt\hbox{$\scriptscriptstyle\square$}\mkern-2mu\raise4pt\hbox{$\scriptscriptstyle\square$}\mkern-2mu\raise8pt\hbox{$\scriptscriptstyle\square$}\mkern1mu}}

\def\OTRUF{\mathinner{\mkern1mu\raise-4pt\vbox{\kern7pt\hbox{$\scriptscriptstyle\square$}}\mkern-2mu\raise0pt\hbox{$\scriptscriptstyle\square$}\mkern-2mu\raise4pt\hbox{$\scriptscriptstyle\square$}\mkern-9.75mu\raise8pt\hbox{$\scriptscriptstyle\square$}\mkern1mu}}

\def\OTRLF{\mathinner{\mkern1mu\raise-4pt\vbox{\kern7pt\hbox{$\scriptscriptstyle\square$}}\mkern-2mu\raise0pt\hbox{$\scriptscriptstyle\square$}\mkern-2mu\raise4pt\hbox{$\scriptscriptstyle\square$}\mkern-17.5mu\raise8pt\hbox{$\scriptscriptstyle\square$}\mkern10.75mu}}

\def\OTURF{\mathinner{\mkern1mu\raise-4pt\vbox{\kern7pt\hbox{$\scriptscriptstyle\square$}}\mkern-2mu\raise0pt\hbox{$\scriptscriptstyle\square$}\mkern-9.75mu\raise4pt\hbox{$\scriptscriptstyle\square$}\mkern-2mu\raise8pt\hbox{$\scriptscriptstyle\square$}\mkern1mu}}

\def\OTURUF{\mathinner{\mkern1mu\raise-4pt\vbox{\kern7pt\hbox{$\scriptscriptstyle\square$}}\mkern-2mu\raise0pt\hbox{$\scriptscriptstyle\square$}\mkern-9.75mu\raise4pt\hbox{$\scriptscriptstyle\square$}\mkern-9.75mu\raise8pt\hbox{$\scriptscriptstyle\square$}\mkern1mu}}

\def\TOTR{\mathinner{\mkern1mu\raise2pt\vbox{\kern7pt\hbox{$\scriptscriptstyle\square$}}\mkern-2mu\raise-2pt\hbox{$\scriptscriptstyle\square$}\mkern-2mu\raise2pt\hbox{$\scriptscriptstyle\square$}\mkern-2mu\raise6pt\hbox{$\scriptscriptstyle\square$}\mkern1mu}}

\def\TOTUR{\mathinner{\mkern1mu\raise2pt\vbox{\kern7pt\hbox{$\scriptscriptstyle\square$}}\mkern-2mu\raise-2pt\hbox{$\scriptscriptstyle\square$}\mkern-2mu\raise2pt\hbox{$\scriptscriptstyle\square$}\mkern-9.75mu\raise6pt\hbox{$\scriptscriptstyle\square$}\mkern1mu}}

\def\RTOf{\mathinner{\mkern1mu\raise8pt\vbox{\kern7pt\hbox{$\scriptscriptstyle\square$}}\mkern-2mu\raise4pt\hbox{$\scriptscriptstyle\square$}\mkern-2mu\raise0pt\hbox{$\scriptscriptstyle\square$}\mkern1mu}}

\def\OTRf{\mathinner{\mkern1mu\raise0pt\vbox{\kern7pt\hbox{$\scriptscriptstyle\square$}}\mkern-2mu\raise4pt\hbox{$\scriptscriptstyle\square$}\mkern-2mu\raise8pt\hbox{$\scriptscriptstyle\square$}\mkern1mu}}

\def\OTURf{\mathinner{\mkern1mu\raise0pt\vbox{\kern7pt\hbox{$\scriptscriptstyle\square$}}\mkern-2mu\raise4pt\hbox{$\scriptscriptstyle\square$}\mkern-9.75mu\raise8pt\hbox{$\scriptscriptstyle\square$}\mkern1mu}}

\def\OTRFf{\mathinner{\mkern1mu\raise0pt\vbox{\kern7pt\hbox{$\scriptscriptstyle\square$}}\mkern-2mu\raise4pt\hbox{$\scriptscriptstyle\square$}\mkern-2mu\raise8pt\hbox{$\scriptscriptstyle\square$}\mkern-2mu\raise12pt\hbox{$\scriptscriptstyle\square$}\mkern1mu}}

\def\OTRUFf{\mathinner{\mkern1mu\raise0pt\vbox{\kern7pt\hbox{$\scriptscriptstyle\square$}}\mkern-2mu\raise4pt\hbox{$\scriptscriptstyle\square$}\mkern-2mu\raise8pt\hbox{$\scriptscriptstyle\square$}\mkern-9.75mu\raise12pt\hbox{$\scriptscriptstyle\square$}\mkern1mu}}

\def\OTRLFf{\mathinner{\mkern1mu\raise0pt\vbox{\kern7pt\hbox{$\scriptscriptstyle\square$}}\mkern-2mu\raise4pt\hbox{$\scriptscriptstyle\square$}\mkern-2mu\raise8pt\hbox{$\scriptscriptstyle\square$}\mkern-17.5mu\raise12pt\hbox{$\scriptscriptstyle\square$}\mkern10.75mu}}

\def\OTURFf{\mathinner{\mkern1mu\raise0pt\vbox{\kern7pt\hbox{$\scriptscriptstyle\square$}}\mkern-2mu\raise4pt\hbox{$\scriptscriptstyle\square$}\mkern-9.75mu\raise8pt\hbox{$\scriptscriptstyle\square$}\mkern-2mu\raise12pt\hbox{$\scriptscriptstyle\square$}\mkern1mu}}

\def\OTURUFf{\mathinner{\mkern1mu\raise0pt\vbox{\kern7pt\hbox{$\scriptscriptstyle\square$}}\mkern-2mu\raise4pt\hbox{$\scriptscriptstyle\square$}\mkern-9.75mu\raise8pt\hbox{$\scriptscriptstyle\square$}\mkern-9.75mu\raise12pt\hbox{$\scriptscriptstyle\square$}\mkern1mu}}

\def\TOTRf{\mathinner{\mkern1mu\raise4pt\vbox{\kern7pt\hbox{$\scriptscriptstyle\square$}}\mkern-2mu\raise0pt\hbox{$\scriptscriptstyle\square$}\mkern-2mu\raise4pt\hbox{$\scriptscriptstyle\square$}\mkern-2mu\raise8pt\hbox{$\scriptscriptstyle\square$}\mkern1mu}}

\def\TOTURf{\mathinner{\mkern1mu\raise4pt\vbox{\kern7pt\hbox{$\scriptscriptstyle\square$}}\mkern-2mu\raise0pt\hbox{$\scriptscriptstyle\square$}\mkern-2mu\raise4pt\hbox{$\scriptscriptstyle\square$}\mkern-9.75mu\raise8pt\hbox{$\scriptscriptstyle\square$}\mkern1mu}}

\def\ODT{\mathinner{\mkern1mu\raise0pt\vbox{\kern7pt\hbox{$\scriptscriptstyle\square$}}\mkern-2mu\raise4.85pt\hbox{$\scriptscriptstyle\times$}\mkern-9.75mu\raise4pt\hbox{$\scriptscriptstyle\square$}\mkern1mu}}

\def\TODT{\mathinner{\mkern1mu\raise4pt\vbox{\kern7pt\hbox{$\scriptscriptstyle\square$}}\mkern-2mu\raise0pt\hbox{$\scriptscriptstyle\square$}\mkern-2mu\raise4.85pt\hbox{$\scriptscriptstyle\times$}\mkern-9.75mu\raise4pt\hbox{$\scriptscriptstyle\square$}\mkern1mu}}

\def\OTDR{\mathinner{\mkern1mu\raise-2pt\vbox{\kern7pt\hbox{$\scriptscriptstyle\square$}}\mkern-2mu\raise2pt\hbox{$\scriptscriptstyle\square$}\mkern-2mu\raise6.85pt\hbox{$\scriptscriptstyle\times$}\mkern-9.75mu\raise6pt\hbox{$\scriptscriptstyle\square$}\mkern1mu}}

\def\ODTR{\mathinner{\mkern1mu\raise-2pt\vbox{\kern7pt\hbox{$\scriptscriptstyle\square$}}\mkern-2mu\raise2.85pt\hbox{$\scriptscriptstyle\times$}\mkern-9.75mu\raise2pt\hbox{$\scriptscriptstyle\square$}\mkern-2mu\raise6pt\hbox{$\scriptscriptstyle\square$}\mkern1mu}}

\def\ODTUR{\mathinner{\mkern1mu\raise-2pt\vbox{\kern7pt\hbox{$\scriptscriptstyle\square$}}\mkern-2mu\raise2.85pt\hbox{$\scriptscriptstyle\times$}\mkern-9.75mu\raise2pt\hbox{$\scriptscriptstyle\square$}\mkern-9.75mu\raise6pt\hbox{$\scriptscriptstyle\square$}\mkern1mu}}

\def\OTRDF{\mathinner{\mkern1mu\raise-4pt\vbox{\kern7pt\hbox{$\scriptscriptstyle\square$}}\mkern-2mu\raise0pt\hbox{$\scriptscriptstyle\square$}\mkern-2mu\raise4pt\hbox{$\scriptscriptstyle\square$}\mkern-2mu\raise8pt\hbox{$\scriptscriptstyle\square$}\mkern-9.75mu\raise8.85pt\hbox{$\scriptscriptstyle\times$}\mkern1mu}}

\def\OTRLDF{\mathinner{\mkern1mu\raise-4pt\vbox{\kern7pt\hbox{$\scriptscriptstyle\square$}}\mkern-2mu\raise0pt\hbox{$\scriptscriptstyle\square$}\mkern-2mu\raise4pt\hbox{$\scriptscriptstyle\square$}\mkern-17.5mu\raise8pt\hbox{$\scriptscriptstyle\square$}\mkern-9.75mu\raise8.85pt\hbox{$\scriptscriptstyle\times$}\mkern10.75mu}}

\def\OTDRF{\mathinner{\mkern1mu\raise-4pt\vbox{\kern7pt\hbox{$\scriptscriptstyle\square$}}\mkern-2mu\raise0pt\hbox{$\scriptscriptstyle\square$}\mkern-2mu\raise4pt\hbox{$\scriptscriptstyle\square$}\mkern-9.75mu\raise4.85pt\hbox{$\scriptscriptstyle\times$}\mkern-2mu\raise8pt\hbox{$\scriptscriptstyle\square$}\mkern1mu}}

\def\OTDRLF{\mathinner{\mkern1mu\raise-4pt\vbox{\kern7pt\hbox{$\scriptscriptstyle\square$}}\mkern-2mu\raise0pt\hbox{$\scriptscriptstyle\square$}\mkern-2mu\raise4pt\hbox{$\scriptscriptstyle\square$}\mkern-9.75mu\raise4.85pt\hbox{$\scriptscriptstyle\times$}\mkern-17.5mu\raise8pt\hbox{$\scriptscriptstyle\square$}\mkern10.75mu}}

\def\ODTURUF{\mathinner{\mkern1mu\raise-4pt\vbox{\kern7pt\hbox{$\scriptscriptstyle\square$}}\mkern-2mu\raise0pt\hbox{$\scriptscriptstyle\square$}\mkern-9.75mu\raise0.85pt\hbox{$\scriptscriptstyle\times$}\mkern-9.75mu\raise4pt\hbox{$\scriptscriptstyle\square$}\mkern-9.75mu\raise8pt\hbox{$\scriptscriptstyle\square$}\mkern1mu}}

\def\ODTURF{\mathinner{\mkern1mu\raise-4pt\vbox{\kern7pt\hbox{$\scriptscriptstyle\square$}}\mkern-2mu\raise0pt\hbox{$\scriptscriptstyle\square$}\mkern-9.75mu\raise0.85pt\hbox{$\scriptscriptstyle\times$}\mkern-9.75mu\raise4pt\hbox{$\scriptscriptstyle\square$}\mkern-2mu\raise8pt\hbox{$\scriptscriptstyle\square$}\mkern1mu}}

\def\ODTRLF{\mathinner{\mkern1mu\raise-4pt\vbox{\kern7pt\hbox{$\scriptscriptstyle\square$}}\mkern-2mu\raise0pt\hbox{$\scriptscriptstyle\square$}\mkern-9.75mu\raise0.85pt\hbox{$\scriptscriptstyle\times$}\mkern-2mu\raise4pt\hbox{$\scriptscriptstyle\square$}\mkern-17.5mu\raise8pt\hbox{$\scriptscriptstyle\square$}\mkern10.75mu}}

\def\ODTRUF{\mathinner{\mkern1mu\raise-4pt\vbox{\kern7pt\hbox{$\scriptscriptstyle\square$}}\mkern-2mu\raise0pt\hbox{$\scriptscriptstyle\square$}\mkern-9.75mu\raise0.85pt\hbox{$\scriptscriptstyle\times$}\mkern-2mu\raise4pt\hbox{$\scriptscriptstyle\square$}\mkern-9.5mu\raise8pt\hbox{$\scriptscriptstyle\square$}\mkern1mu}}

\def\ODTRF{\mathinner{\mkern1mu\raise-4pt\vbox{\kern7pt\hbox{$\scriptscriptstyle\square$}}\mkern-2mu\raise0pt\hbox{$\scriptscriptstyle\square$}\mkern-9.75mu\raise0.85pt\hbox{$\scriptscriptstyle\times$}\mkern-2mu\raise4pt\hbox{$\scriptscriptstyle\square$}\mkern-2mu\raise8pt\hbox{$\scriptscriptstyle\square$}\mkern1mu}}

\def\ODTi{\mathinner{\mkern1mu\raise-2pt\vbox{\kern7pt\hbox{$\scriptscriptstyle\square$}}\mkern-2mu\raise2.85pt\hbox{$\scriptscriptstyle\times$}\mkern-9.75mu\raise2pt\hbox{$\scriptscriptstyle\square$}\mkern-9.75mu\raise6pt\hbox{$\scriptscriptstyle\square$}\mkern-7.5mu\raise12pt\hbox{$\scriptscriptstyle\vdots$}\mkern1mu}}

\def\ODTn{\mathinner{\mkern1mu\raise-10pt\vbox{\kern7pt\hbox{$\scriptscriptstyle\square$}}\mkern-2mu\raise-5.15pt\hbox{$\scriptscriptstyle\times$}\mkern-9.75mu\raise-6pt\hbox{$\scriptscriptstyle\square$}\mkern-9.75mu\raise-2pt\hbox{$\scriptscriptstyle\square$}\mkern-7.5mu\raise4pt\hbox{$\scriptscriptstyle\vdots$}\mkern-7.5mu\raise15pt\hbox{$\scriptscriptstyle\square$}\mkern-1mu\raise6pt\hbox{$\Big\}$}\mkern-1mu\raise7pt\hbox{${\scriptstyle{n}}$}\mkern1mu}}

\def\ODTnm{\mathinner{\mkern1mu\raise-2pt\vbox{\kern7pt\hbox{$\scriptscriptstyle\square$}}\mkern-2mu\raise2.85pt\hbox{$\scriptscriptstyle\times$}\mkern-9.75mu\raise2pt\hbox{$\scriptscriptstyle\square$}\mkern-9.75mu\raise6pt\hbox{$\scriptscriptstyle\square$}\mkern-7.5mu\raise12pt\hbox{$\scriptscriptstyle\vdots$}\mkern-7.5mu\raise23pt\hbox{$\scriptscriptstyle\square$}\mkern-1mu\raise14pt\hbox{$\Big\}$}\mkern-1mu\raise15pt\hbox{${\scriptstyle{n}}$}\mkern1mu}}
\begin{document}


\title{Asymptotic behavior of random heaps}


\author{J. Ben Hough}



\thanks{$^1$Department of Mathematics, UC Berkeley, CA 94720-3860, {\it e-mail:}  jbhough@math.berkeley.edu}





\begin{abstract}
We consider a random walk $W_n$ on the locally free group (or equivalently a signed random heap) with $m$ generators subject to periodic boundary conditions.  Let $\#T(W_n)$ denote the number of removable elements, which determines the heap's growth rate.  We prove that $\lim_{n \rightarrow \infty} \frac{\E(\# T(W_n))}{m} \leq 0.32893$ for $m \geq 4$.  This result disproves a conjecture (due to Vershik, Nechaev and Bikbov [\ref{vershik}]) that the limit tends to $\frac{1}{3}$ as $m \rightarrow \infty$.

\end{abstract}
\maketitle

\newtheorem{theorem}{Theorem}[section]
\newtheorem{corollary}[theorem]{Corollary}
\newtheorem{proposition}[theorem]{Proposition}
\newtheorem{lemma}[theorem]{Lemma}
\newtheorem{defn}[theorem]{Definition}

\newtheorem*{remark}{Remark}
\newenvironment{definition}{\begin{defn}\normalfont}{\end{defn}}

\section{Introduction}
\begin{figure}[htbp]
\centering
\includegraphics[width = 0.25\textwidth]{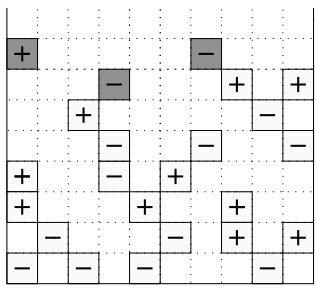}
\caption{A signed random heap with 10 generators.  The roof elements are shaded}
\label{heappic}
\end{figure}

We consider a random walk $W_n$ on the locally free group $\L\F_{m}$ with $m$ generators subject to periodic boundary conditions, $\left< g_0,g_1, \dots, g_{m-1} : g_i g_j = g_j g_i \; \forall \; i \neq j \pm 1 \; {\textrm{mod}}\; m \right>$.  The random walk on this group has an elegant interpretation in terms of random heaps, which are defined in [\ref{vershik1}] and examined more carefully in [\ref{vershik}].  Consider the lattice $\Z_m \times \Z^+$, where $\Z_m$ denotes the integers modulo $m$, and drop signed pieces with $+/-$ charges uniformly over the $m$ columns.  When a piece is dropped in a column it falls as far as possible subject to the condition that it cannot fall past a piece in its own column or a piece in either of the neighboring columns, i.e. ~the pieces have \textquotedblleft sticky corners".   Also, if a piece with a $+$ charge lands directly on top of a piece with a $-$ charge, the two pieces annihilate.  More precisely, each piece in the random heap may be described by three coordinates:  its height, horizontal position and sign.  Let$(h_j, j,s_j)$ denote the highest piece in column $j$, and set $h_j = 0$ if no such piece exists.  When a new piece is to be added to the heap, its horizontal position $k$ and sign $s$ are chosen uniformly over the $2m$ possiblities.  If columns $k-1$, $k$ and $k+1$ are empty, then the height of the new piece is 1.  If $h_k>\max \left\{ h_{k-1},h_{k+1} \right\}$ and $s_k \neq s$, the new piece and $(h_k,k,s_k)$ will both annihilate.  Otherwise, the new piece is added to the heap and its height is $\max \left\{ h_{k-1},h_k,h_{k+1}\right\}+1$.  A typical random heap is illustrated in Figure \ref{heappic}.  

The interpretation of this construction should be clear.  Dropping a $+$ (resp. $-$) piece in the $k^{th}$ column corresponds to adding the generator $g_k$ (resp. $g_k^{-1}$) to the random walk.  The fact that $g_k$ and $g_{k+1}$ don't commute is reflected in the fact that pieces have sticky corners.  In section 2, we describe a bijective correspondence between signed heaps with $m$ generators and reduced words in $\L\F_{m}$ so that a random walk on $\L\F_{m}$ corresponds to the growth procedure of the random heap detailed above.  Our bijection differs from the one presented in [\ref{vershik}], which is for the locally free group without periodic boundary conditions.  This formalism allows us to use the terms random heap process and random walk on $\L\F_{m}$ interchangeably.

This geometrical picture of the random walk on $\L\F_{m}$ leads one to consider the {\it roof} $T(W_n)$ of the random heap, which is the collection of removable pieces.  A piece is said to be {\it removable} if there are no pieces lying above it in its own column or the two adjacent columns.  We shall say that a piece $x$ {\it blocks} piece $y$ if $x$ must be annihilated before $y$ can become part of the roof.  The roof is an important feature of the random heap because its size controls the expected rate of growth of the heap.  Indeed, if the roof has cardinality $k$, then the next piece to fall will annihilate a piece in the heap with probability $\frac{k}{2m}$, and will increase the size of the heap with probability $\frac{2m-k}{2m}$.  The long term rate of growth of the heap is thus controlled by the long term expected size of the roof.  To make this notion precise, we define the {\it drift} of the random heap process to be $\zeta = \lim_{n \rightarrow \infty} \frac{1}{n}\E(\# W_n)$ where $\# W_n$ denotes the number of pieces in the heap.  It follows easily from the above observation that $\zeta = 1 - \lim_{n \rightarrow \infty} \frac{\E(\# T(W_n))}{m}$, a result given in [\ref{vershik}].  

We prove the following:
\begin{theorem}
For $m \geq 4$, $\lim_{n \rightarrow \infty} \frac{1}{m}\E(\# T(W_n)) \leq 0.32893$.
\end{theorem}
This result is surprising, since there are heuristic reasons described in [\ref{vershik}] for suspecting that the limit should converge to $\frac{1}{3}$ as $m \rightarrow \infty$.  In particular, one may consider a random walk on the semi-group $\L\F_{m}^+$ which has generators $g_0, g_1, \dots, g_{m-1}$ satisfying the relations above, but does not include the inverses of the $g_i$'s.  This random walk $\tilde{W}_n$ has an interpretation as an unsigned random heap process.  Specifically, we drop pieces uniformly over $m$ columns as before, but now we forbid pieces from annihilating.  Restricting attention to the roof of the unsigned heap, it is easy to see that this process is an irreducible Markov chain and we have the following:

\begin{proposition}(Proved in [\ref{vershik}])
If $\tilde{W}_n$ is the unsigned random heap process with $m$ generators, then $ \lim_{n \rightarrow \infty}\frac{1}{m} \# T(\tilde{W}_n) = \frac{1}{3}$ a.s.
\end{proposition}

The proof given in [\ref{vershik}] may be simplified.  If one defines $X_{n,i}$ to be the number of roof elements in the $i^{th}$ column after the $n^{th}$ particle is dropped, this process is a Markov chain on $\left\{ 0,1 \right\}$ with transition probabilities $p_{1,0} = 1-p_{1,1}=\frac{2}{m}$ and $p_{0,1} = 1-p_{0,0}=\frac{1}{m}$.  The stationary distribution for this 2-state chain is (2/3,1/3), so $\lim_{n\rightarrow \infty} \E \left( \# T(\tilde{W}_n) \right) = \lim_{n \rightarrow \infty} \sum_{i=1}^m \E(X_{n,i}) = \frac{m}{3}$.  The claim now follows from the ergodic theorem.  This elegant proof was discovered by G\'abor Pete.  From this nice result for unsigned heaps one might suspect that the same stationary behavior should be exhibited by the signed process, at least in the limit $m \rightarrow \infty$, a conjecture expressed in [\ref{vershik}].

The paper is organized as follows.  First we describe the bijection between elements of $\L\F_{m}$ and random heaps.  Then we prove that for $m \geq 4$, 
\begin{equation}
	\limsup_{n \rightarrow \infty} \frac{1}{mn} \sum_{k = 1}^n \# T(W_k) \leq 0.32893 \; a.s. \label{time}
\end{equation}
Finally we shall deduce that if one allows heaps to be infinite, the random heap process has a unique stationary distribution $\nu$ and the finite dimensional distributions of $W_n$ converge to the f.d.d.'s of $\nu$ as $n \rightarrow \infty$.  From this fact, we deduce the same upper bound for the space average, proving the theorem.  We conclude with some open problems. 

\section{The Bijection between Heaps and Words}
We begin by giving a formal definition of a heap, and introducing some notation.
\begin{definition} \label{heapdef}
A {\it heap} $H$ is a finite union of pieces $(v_i,w_i,\sigma_i)$ with $v_i \in \Z^+$, $w_i \in \Z_m$ and $\sigma_i \in \left\{+,-\right\}$, which satisfy the following conditions:
\begin{itemize}
\item[1.]  If $(v_\alpha,w_\alpha,\sigma_\alpha) \in H$ and $v_\alpha>1$, then there exists $(v_\beta,w_\beta,\sigma_\beta)\in H$ such that $v_\beta= v_\alpha-1$ and $w_\beta\in \left\{w_\alpha-1,w_\alpha,w_\alpha+1\right\}$.
\item[2.]  If $(v_\alpha,w_\alpha,\sigma_\alpha), (v_\beta,w_\beta,\sigma_\beta)\in H$ satisfy $v_\alpha= v_\beta+1$ and $w_\alpha= w_\beta$, then $\sigma_\alpha= \sigma_\beta$.
\item[3.]  If $(v_\alpha,w_\alpha,\sigma_\alpha), (v_\beta,w_\beta,\sigma_\beta)\in H$ satisfy $w_\alpha= w_\beta+1$, then $v_\alpha\neq v_\beta$.
\end{itemize} 
\end{definition}
If $H$ is a heap, the elements of the roof may all be removed to construct a new heap $H'$.  We shall say that the roof of $H'$ is the $2^{nd}$ level roof of $H$.  Similarly, we may define the $3^{rd}$ level roof and so forth.  It is easy to see that a heap may be determined uniquely by specifying the horizontal positions of the pieces in each of its roofs.

The concept of a roof also makes sense for words representing elements of $\L\F_{m}$.  Specifically, if $a$ is an instance of $g_i$ or $g_i^{-1}$ in the word $W$, we say that $a$ is  {\it removeable} if it commutes with all the letters in $W$ occuring to its right, and no other instances of $g_i$ or $g_i^{-1}$ occur to its right.  The {\it roof} of $W$ consists of the collection of removeable letters.  Omitting the removeable letters from $W$ to obtain $W'$, we define the roof of $W'$ to be the $2^{nd}$ level roof of $W$.  Higher level roofs are defined analogously.  If $W = e$, the unit in $\L\F_{m}$, we say that all the roofs are empty.  Observe that two words with the same roofs must represent the same element in $\L\F_{m}$, although words with different roofs may also be equal.  To obtain a unique representation, we introduce the following definition.

\begin{figure}[htbp]
\centering
\includegraphics[width = 0.4\textwidth]{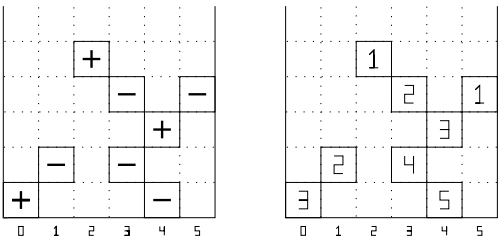}
\caption{The pieces in the second figure are numbered to indicate their roof membership.  The word corresponding to this heap is $g_4^{-1}g_3^{-1} g_0 g_4 g_1 g_3 g_2 g_5$}
\label{bijpic}
\end{figure}

\begin{definition}
A word $W$ representing an element of $\L\F_{m}$ is said to be in {\it normal form} if the following conditions are satisfied.  Here $a$ is an instance of $g_i$ or $g_i^{-1}$ in $W$, and $b$ is an instance of $g_j$ or $g_j^{-1}$.    
\begin{itemize}
\item[1.] If $ab = e$, then $a$ and $b$ are not members of the same or adjacent roofs.
\item[2.] If $a$ and $b$ are members of the same roof and $i<j$, then $b$ occurs to the right of $a$ in $W$.
\item[3.] If $a$ is a member of a higher level roof than $b$, then $b$ occurs to the right of $a$ in $W$.
\end{itemize}
\end{definition}
The fact that every element of $\L\F_{m}$ may be represented by a unique normalized word follows directly from the commutation relations defining $\L\F_{m}$.  By identifying the respective roofs of a heap and a normalized word so that a $+$ (resp. $-$) piece in column $k$ corresponds to an instance of the generator $g_k$ (resp. $g_k^{-1}$) one obtains the desired bijection, see Figure \ref{bijpic}.  With this identification, the random growth process of a heap described in the introduction corresponds to a random walk on $\L\F_{m}$.

\section{Bounding the Time Average}
We now give a proof of inequality (\ref{time}).  Consider a single column of a random heap, and when each piece is dropped record a \textquotedblleft 1" if the column contains an element of the roof and a \textquotedblleft 0" otherwise.  This generates a sequence of 0's and 1's, say $X_1, X_2, X_3, \dots$, and to prove (\ref{time}) it is enough to show that $\limsup_{n \rightarrow \infty} \frac{1}{n} \sum_{k=1}^n X_k \leq 0.32893 \; {\textrm{a.s.}}$  The sequence $X_1, X_2, X_3 \dots$ may be expressed more compactly as $\0_0, \1_1, \0_1, \1_2, \0_2, \dots$ where $\1_k$ denotes the length of the $k^{th}$ sequence of 1's and similarly for $\0_k$.  The random variable $\0_0$ is distinguished in that its value may be zero.  Our approach is to show that 
\begin{equation}
\limsup_{n \rightarrow \infty} \frac{\sum_{k=1}^n \1_k}{\sum_{k=1}^n \1_k + \sum_{k=1}^{n-1} \0_k} \leq 0.32893 \; {\textrm{a.s.}} \label{limsup}
\end{equation}
We shall do this by constructing ${\textrm{a.s.}}$ upper and lower bounds for 
\begin{equation}
\limsup_{n \rightarrow \infty} \frac{1}{n} \sum_{k=1}^\infty \1_k \;\;  {\textrm{and}}\;\; \liminf_{n\rightarrow \infty} \frac{1}{n} \sum_{k=1}^n \0_k,
\end{equation}
 respectively.

An observation fundamental to our proof is that a sequence of 0's may be terminated in two different ways.
\begin{itemize}
\item[1.] A piece may fall in the distinguished column.
\item[2.] All the pieces blocking the highest piece in the distinguished column may be annihilated.  
\end{itemize}
We shall say that a 0 sequence {\it builds upward} if it is terminated via method 1 and otherwise it {\it backtracks}.  Similarly, a sequence of 1's may be terminated in two ways.
\begin{itemize}
\item[1.]  A piece may fall in a column adjacent to the distinguished column.
\item[2.]  All the pieces in the distinguished column blocking the highest piece in an adjacent column may be annihilated (or all the pieces in the distinguished column may be annihilated if the adjacent columns contain no pieces).
\end{itemize}
As before, we say that a sequence of 1's builds upward if it is terminated by method 1, and that it backtracks otherwise.

The first task is to construct an a.s. ~upper bound for $\limsup_{n \rightarrow \infty} \frac{1}{n} \sum_{k=1}^n \1_k$, the following outlines our approach.  The key idea is based on the observation above.  In general, if a sequence of 0's builds upward we expect the subsequent sequence of 1's to be significantly shorter than if the sequence of 0's backtracks.  For in the first case the subsequent sequence of 1's may be easily terminated via either backtracking or building upward.  Whereas in the latter case the sequence of 1's will likely begin at a time when there is a stack of pieces directly below the roof element in the distinguished column so it will be much more difficult to terminate the sequence of 1's via backtracking.  We will describe this phenomenon precisely in Lemma \ref{lemma1}, then will show in Lemma \ref{plemma} that 0 sequences usually build upward.  Henceforth, we shall say that a sequence of 1's starts from a {\it short position} if it is initiated at a time when there is no piece directly below the roof piece in the distinguished column.  Otherwise we say that the 1 sequence starts from a {\it long position}.  Also, we shall refer to the distinguished column as column 0 and label the columns to the right $1,2,3, \dots$ modulo $m$.  Unless otherwise stated, we assume without comment that $m \geq 4$.

\begin{lemma}  
$\E(\1_k | \1_k \textrm{ starts from a short position}) = (\sqrt{2}-1)m$ and $\E(\1_k) \leq \frac{m}{2}$. \label{lemma1}
\end{lemma}
\begin{proof}
Consider the following random walk procedure on $\Z$.  Start the walk from $\ell>0$, and between consecutive steps of the walk wait ${\textrm{i.i.d.}}$ times $\tau_i$ with distribution $\P(\tau_i>s) = \left(\frac{m-3}{m}\right)^s$.  Also, immediately prior to taking each step flip a coin and with probability 2/3 stop the walk.  If the walk is not stopped, then move 1 unit to the left or right with equal probability.  Stop the walk when it reaches the origin, if it has not been stopped previously.  It is easy to see that this process is equivalent to the process that determines the length of $\1_k$.  The position of the random walk corresponds to the number of consecutive pieces at the top of column 0.  Setting $\ell = 1$ we find:
\begin{equation}
\E(\1_k | \1_k \textrm{ starts from a short position}) = \E(\tau_i)\E(N \wedge T),
\end{equation}
where $\P(N>k) = \left(\frac{1}{3}\right)^k$ and $T$ is the hitting time of 0 for a simple random walk starting from 1.  It follows from Lemma 3.3 in Chapter 3 of [\ref{durrett}] that
\begin{equation}\
\P(T > 2k-1) = \P(T > 2k) = {2k \choose k} \left( \frac{1}{4}\right)^k 
\end{equation}
so, since $N$ and $T$ are independent,
\begin{equation}
\E(N \wedge T) = 1+ \sum_{k=1}^\infty {2k \choose k} \left(\frac{1}{4}\right)^k \left[\left(\frac{1}{3}\right)^{2k-1} + \left(\frac{1}{3}\right)^{2k} \right] = 3(\sqrt{2}-1).
\end{equation}
The above sum is evaluated using the identity $\sum_{k=0}^\infty {2k \choose k} t^k =(1-4t)^{-1/2}$ which is valid for $|t|<\frac{1}{4}$.  Since $\E(\tau_i) = \frac{m}{3}$ we obtain $\E(\1_k) = (\sqrt{2}-1)m$.  From any starting position, $\E(\1_k) \leq \E(\tau_i)\E(N) = \frac{m}{2}$.
\end{proof}

\begin{figure}[htbp]
\centering
\includegraphics[width = 0.4\textwidth]{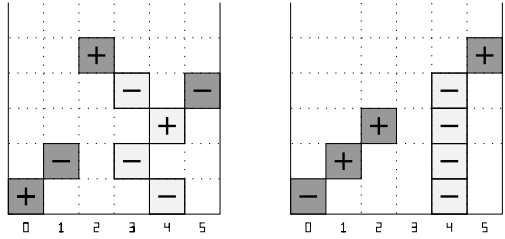}
\caption{The shaded pieces in each figure are both realizations of the symbol $\TOTR$}
\label{heappic2}
\end{figure}

Next we bound the probability that a sequence of 0's backtracks.  We say that a sequence of 0's starts from a short position if it is initiated at a time when there is only one piece blocking the highest piece in column zero.  Then it is clear that the probability of backtracking is maximized when the 0 series starts from a short position.  In what follows, we write 
\begin{equation}
 	\OT, \TO, \TOT, \RTOf, \OTRf,  \OTURf,   \TOTRf, \TOTURf, \OTRFf, \OTRUFf, \OTRLFf, \OTURFf, \OTURUFf
\end{equation}
to denote possible configurations for the pieces blocking the uppermost piece in column zero.  In each of the above pictures, the lowest piece is assumed to lie in column zero and we agree that exactly those pieces which block this piece are shown.  The position of a piece indicates its column and which of the other pieces it blocks.  It does not necessarily indicate its height relative to the others, see Figure \ref{heappic2}.  Two different realizations of a given symbol are said to represent the same {\it configuration} of blocking pieces.  We shall use the notation $\P(\TOT)$ to denote the probability that a 0 sequence backtracks from the starting position $\TOT$, and so forth.  To see that this notation makes sense, we appeal to the following lemma.

\begin{lemma}
The configuration of pieces blocking the highest piece in column zero determines both the probability of backtracking and the probability that this piece will be annihilated at any later time. \label{basiclemma}
\end{lemma}
\begin{proof}
Consider random heaps $W^1$ and $W^2$, and suppose that at time $n$ the pieces blocking the highest piece in column zero of $W^1_n$ and $W^2_n$ have the same configuration.  Label the highest pieces in column zero $x_1$ and $x_2$ respectively.  Construct a coupling between $W^1$ and $W^2$ as follows.  Select a piece uniformly over the $2m$ possible choices and add it to $W^1_n$ to obtain $W^1_{n+1}$.  Now add a piece in the same column to $W^2_n$, and choose its sign according to the following rules.
\begin{itemize}
\item[1.]  If the newly added pieces both land on top of roof elements, fix the sign of the second piece so that both annihilate or neither annihilate.
\item[2.]  Otherwise, choose the sign of the second piece to be the same as the sign of the first piece.
\end{itemize}
Observe that at time $n+1$, either both $x_1$ and $x_2$ will have been annihilated or the pieces blocking $x_1$ and $x_2$ will again have the same configuration.  Iterating this procedure, the proof is complete.
\end{proof}

We now estimate $\P(\OT)$.  Observe that by symmetry we have $\P(\OT) = \P(\TO)$.

\begin{lemma}
$0.137457 \leq \P(\OT) \leq 0.14599$. \label{plemma}
\end{lemma}
\begin{proof}
Conditioning on the next piece to fall in either column -1,0,1 or 2 we obtain:
\begin{equation}
\P(\OT) = \frac{1}{4}\P(\TOT)+\frac{1}{4}\cdot 0 + \frac{1}{8}\cdot 1 + \frac{1}{8}\P\left(\OTUR\right) + \frac{1}{4} \P\left(\OTR\right).
\end{equation}
Now expanding the terms on the RHS and using symmetry:
\begin{eqnarray}
\frac{291}{320}\P(\OT) &=& \frac{1}{8} + \frac{13}{160} \P\left(\TOTUR\right) + \frac{3}{20}\P\left(\TOTR\right) + \frac{1}{20}\P\left(\OTRF\right) + \frac{1}{20}\P\left(\OTRLF\right)  \label{Pexp}\\
& & + \frac{1}{64} \P\left(\OTURUF\right) + \frac{1}{32} \P\left(\OTURF\right) + \frac{1}{40}\P\left(\OTRUF\right).\nonumber
\end{eqnarray}
The lower bound for $\P(\OT)$ is obtained by setting the unknown $\P(\cdot)$ terms on the RHS equal to zero.  For the upper bound, we use the following lemma.
\begin{lemma}
For any initial heap, call $\wp$ the probability that a specific piece in the roof is annihilated at any future time.  Then $\frac{1}{6}<\wp<\frac{1}{5}$. \label{blemma}
\end{lemma}
\begin{proof}
Conditioning on the next piece to fall in any of the three relevant columns, we obtain
\begin{equation}
\frac{1}{6}<\wp < \frac{1}{6} + \frac{5}{6}\wp^2.
\end{equation}
Thus $(5\wp - 1)(\wp-1)>0$, and we deduce that $\wp<\frac{1}{5}$.  
\end{proof}
It follows that 
\begin{eqnarray}
\P\left(\OTRF\right), \P\left(\OTURUF\right), \P\left(\OTRLF\right),\P\left(\OTURF\right) &\leq& \frac{1}{25} \P(\OT), \nonumber \\
\P\left(\TOTR\right), \P\left(\TOTUR\right) &\leq& \frac{1}{5} \P(\OT).
\end{eqnarray}
Substituting these approximations into (\ref{Pexp}) yields the upper bound.
\end{proof}

The next task is to construct an ${\textrm{a.s.}}$ upper bound on $\limsup_{n\rightarrow \infty}\frac{1}{n} \sum_{k=1}^n \1_k$ where $\0_0,\1_1,\0_1,\1_2, \dots$ are chosen according to the heap process.  We construct $\tilde{\1}_1, \tilde{\1}_2, \tilde{\1}_3, \dots$ so that $\tilde{\1}_k \geq \1_k$ for all $k$ and $\limsup_{n \rightarrow \infty}\frac{1}{n} \sum_{k=1}^n \tilde{\1}_k$ has a sufficiently small ${\textrm{a.s.}}$ upper bound.  The $\tilde{\1}_k$'s are constructed as follows.  Consider the following four random variables: $S$, $S^*$, $\tilde{S}$ and $L$.  The variable $S$ shall be given the distribution of $\1_k$ starting from a short position.  We give $S^*$ the distribution of $S$ conditioned on the event that the first piece to fall in either column -1, 0 or 1 after the 1 sequence is initiated actually falls in either column -1 or 1.  Also, $\tilde{S}$ has the distribution of $S$ conditioned on the event that the first piece to fall in either column -1, 0 or 1 after the 1 sequence is initiated actually falls in column 0.  Finally, $L$ has the distribution of a $\1_k$ whose starting position consists of an infinite tower of pieces in column 0.  From Lemma \ref{lemma1}, $\E(S) = (\sqrt{2} - 1)m$ and $\E(S^*) = \E(\tau_i) = \frac{m}{3}$, moreover $\E(L) = \E(N) = \frac{m}{2}$.  Since $\frac{1}{3} \E(\tilde{S}) + \frac{2}{3} \E(S^*) = \E(S)$, we obtain $\E(\tilde{S}) = (3\sqrt{2} - 11/3)m$.

Suppose now that we have an infinite number of $\mathrm{i.i.d.}$ copies of $S^*$, $\tilde{S}$ and $L$.  To construct $\tilde{\1}_k$, we first build a Markov chain $\xi_k$ on the state space $\left\{s^*,\tilde{s},\ell\right\}$, and then define $\tilde{\1}_k$ to be a new independent copy of $S^*$, $\tilde{S}$ or $L$ if $\xi_k$ equals $s^*$, $\tilde{s}$ or $\ell$ respectively.  We shall ignore the condition $\tilde{\1}_k \geq \1_k$ for the time being.  To construct $\xi_k$, take $\xi_1 = s^*$ with probability 2/3 and $\xi_1 = \tilde{s}$ with probability 1/3.  Once $\xi_k$ has been chosen, choose $\xi_{k+1}$ according to the following transition probabilities.  Here $\rho = \P(\OT)$.
\begin{equation}
s^*\begin{array}{ccc}
  & s^* & \frac{2}{3} \\
\nearrow&  &  \\
\rightarrow & \tilde{s} & \frac{1}{3} \\
\searrow  &  &  \\
  & \ell & 0
\end{array}
\;\;\;\;\;\;\;\;\;\;\;\;\;\;\;\;\;\;
\tilde{s}\begin{array}{ccc}
  & s^* & \frac{2}{3}(1-\rho) \\
\nearrow&  &  \\
\rightarrow & \tilde{s} & \frac{1}{3}(1-\rho) \\
\searrow  &  &  \\
  & \ell & \rho
\end{array}
\;\;\;\;\;\;\;\;\;\;\;\;\;\;\;\;\;\;
\ell\begin{array}{ccc}
  & s^* & \frac{2}{3}(1-\rho) \\
\nearrow&  &  \\
\rightarrow & \tilde{s} & \frac{1}{3}(1-\rho) \\
\searrow  &  &  \\
  & \ell & \rho
\end{array}
\end{equation}
Then $\xi_k$ has stationary distribution:
\begin{equation}
\pi_{s^*} = \frac{2(1-\rho)}{3-2\rho}, \;\;\;\;\;\;\; \pi_{\tilde{s}} = \frac{1-\rho}{3-2\rho}, \;\;\;\;\;\;\; \pi_\ell = \frac{\rho}{3-2\rho}. \label{stat}
\end{equation}
Since all the copies of $S^*$ appearing among the $\tilde{\1}_k$'s are ${\textrm{i.i.d.}}$, and similarly for the copies of $\tilde{S}$ and $L$, it follows from the ergodic theorem that
\begin{equation}
	\frac{1}{n} \sum_{k=1}^n \tilde{\1}_k \ascv \pi_{s^*} \E(S^*) + \pi_{\tilde{s}}\E(\tilde{S}) + \pi_\ell \E(L).\label{ublimsup}
\end{equation}
Combining equations (\ref{stat}) and (\ref{ublimsup}) with Lemma \ref{plemma} we obtain
\begin{equation}
	\limsup_{n \rightarrow \infty} \frac{1}{n} \sum_{k=1}^n \tilde{\1}_k \leq 0.41884m \;{\textrm{a.s.}}
\end{equation}

Now, to deduce the same ${\textrm{a.s.}}$ bound for $\limsup_{n \rightarrow \infty} \frac{1}{n} \sum_{k=1}^n \1_k$, it remains only to demonstrate that the $\tilde{\1}_k$'s can be constructed so that $\tilde{\1}_k \geq \1_k$.  
\begin{lemma}
	We may construct the sequences $\xi_k$ and $\tilde{\1}_k$ with the joint densities specified above so that $\1_k \leq \tilde{\1}_k$, and if $\xi_k = s^*$, then after $\1_k$ is initiated a piece falls in either column 1 or -1 before one falls in column 0.
\end{lemma}
\begin{proof}
The proof is by induction.  For the base case, take $\tilde{\1}_1 = \1_1$.  To determine $\xi_1$, consider the first piece to fall in either columns -1, 0 or 1 after $\1_1$ is initiated.  If it falls in either column -1 or 1, then $\xi_1 = s^*$, otherwise $\xi_1 = \tilde{s}$.  Now assume that $\xi_1$, \dots, $\xi_n$ and $\tilde{\1}_1, \dots, \tilde{\1}_n$ have been constructed to satisfy the stipulations specified above.  If $\xi_n = s^*$ we see by the induction hypothesis that $\1_{n+1}$ must start from a short position regardless of whether $\0_n$ backtracks.  So take $\tilde{\1}_{n+1} = \1_{n+1}$ in this case.  Define $\xi_{n+1}$ as before by considering the first piece to fall in either column -1, 0 or 1 after $\1_{n+1}$ is initiated.  

If $\xi_n = \tilde{s}$ or $\xi_n = \ell$, then construct a random variable $\chi_{n+1}$ independent of $\tilde{\1}_1, \dots, \tilde{\1}_n$ and $\xi_1, \dots, \xi_n$ so that $\P(\chi_{n+1} = 1) = 1-\rho$ and $\P(\chi_{n+1} = 0) = \rho$.  If $\H_n$ denotes the configuration of the heap at the time when $\0_n$ is initiated, we know from the discussion immediately preceding Lemma \ref{plemma} that $\P(\0_n \textrm{ backtracks }|\H_n) \leq \rho$.  Thus, we may construct $\chi_{n+1}$ so that $\left\{\0_n \textrm{ backtracks}\right\} \subset \left\{ \chi_{n+1} = 0 \right\}$.  Now, if $\chi_{n+1} = 1$, take $\tilde{\1}_{n+1} = \1_{n+1}$.  Construct $\xi_{n+1}$ as before by considering the first piece to fall in either column -1, 0 or 1 after $\1_{n+1}$ is initiated.  If $\xi_n = 0$, take $\xi_{n+1} = \ell$.  To construct $\tilde{\1}_{n+1}$, start building a random heap process on top of an infinite stack of pieces in the 0 column.  Drop pieces in the same order as they fell while $\1_{n+1}$ was being constructed until $\1_{n+1}$ was terminated.  At this point, if the 0 column still contains a piece in the roof, continue to drop signed pieces uniformly over the $m$ columns, but now independently of the random heap process.  Take $\tilde{\1}_{n+1}$ to be the length of the sequence of 1's that is generated.  It is clear that if $\1_{n+1}$ builds upward, it will equal $\tilde{\1}_{n+1}$, and otherwise it will be strictly smaller.
\end{proof}
Thus, we have established the following estimate:
\begin{lemma}
	$\limsup_{n \rightarrow \infty} \frac{1}{n} \sum_{k=1}^n \1_k \leq 0.41884m$ a.s.
\end{lemma}

We now calculate an a.s. lower bound for $\liminf_{n \rightarrow \infty} \frac{1}{n} \sum_{k=1}^n \0_k$.  The first step is to estimate the conditional expectation $\E(\0_k | \0_k \textrm{ starts from a short position})$.  In what follows, we shall write $\E(\TOT)$ to denote $\E(\0_k | \0_k \textrm{ starts from} \TOT)$.  
\begin{lemma}
$0.85453m \leq \E(\OT) \leq 0.86255m$. \label{lemmaEb}
\end{lemma}
\begin{proof}
Conditioning on the first piece to fall in either column -1,0,1 or 2 we obtain
\begin{equation}
	\E(\OT) = \frac{m}{4} + \frac{1}{4}\E(\TOT) + \frac{1}{8}\E\left(\OTUR\right) + \frac{1}{4}\E\left(\OTR\right),
\end{equation}
and expanding the terms on the RHS again yields:
\begin{eqnarray}
\frac{291}{320}\E(\OT) &=& \frac{61m}{160} + \frac{1}{64} \E\left(\OTURUF\right) + \frac{13}{160}\E\left(\TOTUR\right) + \frac{1}{32} \E\left(\OTURF\right) + \frac{3}{20}\E\left(\TOTR\right)  \label{E1}\\
& & + \frac{1}{20}\E\left(\OTRLF\right) + \frac{1}{20}\E\left(\OTRF\right) + \frac{1}{40}\E\left(\OTRUF\right). \nonumber
\end{eqnarray}
Bounding the unknown terms on the RHS by $m$ yields the upper bound.

For the lower bound, it is useful to first estimate the related quantities $\tilde{\P}(\ODT)$ and $\tilde{\E}(\ODT)$.  Here, $\tilde{\P}(\ODT)$ is the probability that starting from the indicated configuration, the piece distinguished by the $\times$ symbol will be annihilated before a piece falls in column 0.  Similarly, $\tilde{\E}(\ODT)$ gives the expected number of steps before either the distinguished piece is annihilated or a piece falls in column 0.  
\begin{lemma}
We have the estimates $0.180115 \leq \tilde{\P}(\ODT) \leq 0.1806355$ and $0.133939 \leq \tilde{\P}\left(\OTDR\right)\leq 0.141677$.  For $m = 4$, $\tilde{\P}\left(\OTRDF\right) = \tilde{\P}(\ODT)$ and for $m \geq 5$ we have $0.133939 \leq \tilde{\P}\left(\OTRDF\right) \leq 0.141677$.  \label{lemmap}
\end{lemma}
\begin{proof}
Conditioning on the next piece to fall in column -1, 0, 1 or 2 we write:
\begin{eqnarray}
	\tilde{\P}(\ODT) &=& \frac{1}{4}\tilde{\P}(\TODT)+ \frac{1}{4} \cdot 0 + \frac{1}{8}\cdot 1 + \frac{1}{8}\tilde{\P}\left(\ODTUR\right)  + \frac{1}{4}\tilde{\P}\left(\ODTR\right) \nonumber  \\
&=& \frac{1}{4}\tilde{\P}(\ODT)+\frac{1}{8} + \frac{1}{8}\tilde{\P}(\ODT)^2  + \frac{1}{4}\tilde{\P}\left(\OTDR\right)\tilde{\P}(\ODT). \label{pt0}
\end{eqnarray}
Expanding $\tilde{\P}\left(\OTDR\right)$ by conditioning on the next piece to fall in column 0,1,2 or 3 gives:
\begin{equation}
\tilde{\P}\left(\OTDR\right) = \frac{1}{4} \tilde{\P}\left(\OTDRLF\right)+\frac{1}{8} + \frac{1}{8}\tilde{\P}\left(\OTDR\right)^2  + \frac{1}{4}\tilde{\P}\left(\OTDRF\right). \label{pt1}
\end{equation}
The probability that the highest piece in $\OTDRF$ is annihilated at all is $\wp \leq 1/5$, and the probability that it is annihilated by the first piece to fall in columns 0, 1, 2, 3 or 4 is 1/10 (actually with 1/10 replaced by 1/8 if $m=4$, but this extra precision is not important).  Thus, we obtain the inequality:
\begin{equation}
\frac{1}{10}\tilde{\P}\left(\OTDR\right) \leq \tilde{\P}\left(\OTDRF\right) \leq \frac{1}{5} \tilde{\P}\left(\OTDR\right). \label{pt3}
\end{equation}
The quantity $\P\left(\OTDRLF\right)$ satisfies the same inequality.  Combining these inequalities with (\ref{pt1}) and (\ref{pt3}), the bounds on $\tilde{\P}\left(\OTDR\right)$ follow.  Once these inequalities are established, the bounds on $\tilde{\P}(\ODT)$ follow from (\ref{pt0}).  To bound $\tilde{\P}\left(\OTRDF\right)$ if $m \geq 5$, condition on the next piece to fall in column 0, 2, 3 or 4.  Then use the analogues of (\ref{pt3}) to bound the unknown higher order terms.
\end{proof}

The quantity $\tilde{\E}(\ODT)$ may be readily calculated from $\tilde{\P}(\ODT)$.  Indeed, let $T$ be the number of steps until either the distinguished piece is annihilated or a piece falls in column zero.  By applying the strong Markov property at $T$, we obtain
\begin{equation}
	\tilde{\E}\left( \ODTUR\right) = \left[1 + \tilde{\P}(\ODT) \right] \tilde{\E}(\ODT) \label{Efact1}.
\end{equation}
Iterating this procedure yields
\begin{equation}
\tilde{\E}\left(\ODTn\right) = \left[1 + \tilde{\P}(\ODT) + \dots + \tilde{\P}(\ODT)^n \right] \tilde{\E}(\ODT),
\end{equation}
and letting $n \rightarrow \infty$ we deduce that 
\begin{equation}
\tilde{\E}(\ODT) = \left(1-\tilde{\P}(\ODT) \right)m \label{Efact2}.
\end{equation}
Analogous expressions clearly hold for $\tilde{\E}\left(\OTDR\right)$ and $\tilde{\E}\left(\OTRDF\right)$.  Using (\ref{Efact2}) and conditionings similar to the one leading to (\ref{Efact1}), we are now able to give lower bounds for the unknown terms in (\ref{E1}).  
\begin{eqnarray}
\E\left(\OTURUF\right) &\geq& \tilde{\E}\left(\ODTURUF\right) = \left[1 +  \tilde{\P}(\ODT)+  \tilde{\P}(\ODT)^2 \right]\tilde{\E}(\ODT) \geq 0.994106m \\
\E(\TOTUR)&\geq& \tilde{\E}\left(\ODTUR\right) \geq \left[1 +  \tilde{\P}(\ODT)\right]\tilde{\E}(\ODT) \geq 0.96737m \\
\E(\OTURF)&\geq& \tilde{\E}\left(\ODTURF\right) \geq \tilde{\E}(\OTDR) + \tilde{\P}\left(\OTDR\right)\left[1 +  \tilde{\P}(\ODT)\right]\tilde{\E}(\ODT) \geq 0.995377m \\ 
\E(\TOTR) &\geq& \tilde{\E}\left(\ODTR\right) \geq \tilde{\E}\left(\OTDR\right) + \tilde{\P}\left(\OTDR\right)\tilde{\E}(\ODT) \geq 0.974408m\\
\E(\OTRLF) &\geq& \tilde{\E}\left(\ODTRLF\right) \geq \tilde{\E}(\ODT) + \tilde{\P}(\ODT)\left[\tilde{\E}\left(\OTDR\right) + \tilde{\P}\left(\OTDR\right)\tilde{\E}(\ODT)\right]\geq 0.995377m \\
\E\left(\OTRF\right) &\geq& \tilde{\E}\left(\ODTRF\right) \geq \tilde{\E}(\OTRDF) + \tilde{\P}(\OTRDF)\left[\tilde{\E}\left(\OTDR\right) + \tilde{\P}\left(\OTDR\right)\tilde{\E}(\ODT)\right]\geq 0.995377m \\
\E\left(\OTRUF\right) &\geq& \tilde{\E}\left(\ODTRUF\right) \geq \tilde{\E}(\OTDR) + \tilde{\P}(\OTDR)\left[\tilde{\E}\left(\OTDR\right) + \tilde{\P}\left(\OTDR\right)\tilde{\E}(\ODT)\right]\geq 0.996374m 
\end{eqnarray}
The numerical estimates are generated by using the estimates from Lemma \ref{lemmap}.  Substituting these values into equation (\ref{E1}) now gives the lower bound for $\E(\OT)$.
\end{proof}

It remains to construct an a.s. ~lower bound for $\liminf_{n \rightarrow \infty} \frac{1}{n} \sum_{k=1}^n \0_k$.  Let $\0$ be a random variable with distribution equal to $\0_k$ starting from a short position.  Since the distribution of $\0_k$ given $\H_k$ is stochastically greater than $\0$, we can construct $\tilde{\0}_1, \tilde{\0}_2, \tilde{\0}_3, \dots$ to be ${\textrm{i.i.d.}}$ copies of $\0$ so that $\tilde{\0}_k \leq \0_k$ for all $k$.  By Lemma \ref{lemmaEb} and the strong law of large numbers, we know that $\limsup_{n\rightarrow \infty}\frac{1}{n} \sum_{k=1}^n \tilde{\0}_k \geq 0.85453m$ a.s.  Hence, we arrive at the estimate
\begin{lemma}
$\limsup_{n\rightarrow \infty} \frac{1}{n} \sum_{k=1}^n \0_k \geq 0.85453m \; \textrm{a.s.}$
\end{lemma}
It follows that 
\begin{eqnarray}
	\limsup_{n \rightarrow \infty} \frac{1}{n} \sum_{k=1}^n X_k &\leq& \limsup_{n \rightarrow \infty} \frac{\sum_{k=1}^n \1_k}{\sum_{k=1}^n \1_k + \sum_{k=1}^{n-1} \0_k} \\
&\leq& \frac{0.41884}{0.41884 + 0.85459} \leq 0.32893. \nonumber 
\end{eqnarray}
By symmetry, this result is valid for any column, and equation (\ref{time}) follows.

\section{Bounding the Space Average}
We have shown that $\limsup_{n \rightarrow \infty} \frac{1}{mn} \sum_{k=1}^n \#T(W_k) \leq 0.32893$ a.s.  Since $0 \leq \#T(W_k) \leq m$ for all $k$, it follows that
\begin{equation}
\limsup_{n\rightarrow \infty} \frac{1}{mn} \sum_{k=1}^n \E(\#T(W_k)) \leq 0.32893.
\end{equation}
Thus, to prove that $\limsup_{n \rightarrow \infty} \frac{1}{m}\E(\#T(W_n)) \leq 0.32893$ it suffices to show that $\lim_{n\rightarrow \infty} \E(\#T(W_n))$ exists.  We shall prove the much stronger result that if one allows heaps to be infinite, the random heap process has a unique stationary distribution $\nu$ and the finite dimensional distributions of $W_n$ converge to the corresponding f.d.d.'s of $\nu$ in total variation distance as $n \rightarrow \infty$.

We begin by introducing some notation.  We shall write $T_k(W)$ to denote the configuration of the roofs of $W$ above the $(k+1)^{st}$ level.  Thus, $T_1(W) = T(W)$, $T_2(W)$ consists of the $1^{st}$ and $2^{nd}$ level roofs, etc.  Infinite heaps are determined abstractly by specifying the configuration of their roof at each level.  The roofs of an infinite heap $W_\infty$ must satisfy the same geometrical relations as those of a finite heap:  $T_k(W_\infty)$ must determine a finite heap satisfying Definition \ref{heapdef} for all $k$.  Finite heaps are a subset of these generalized heaps distinguished by the property that only finitely many of their roofs are nonempty.  Since the uppermost $k$ roofs of a heap are sufficient to determine the uppermost $k-1$ roofs after a piece is added, the random heap process extends naturally to a Markov operator, $\PP$, on the space $\Omega_H$ of generalized heaps.   In what follows $R_k$ shall denote the space of possible configurations of the uppermost $k$ roofs, and $\MM(R_k)$ shall denote the metric space of probability measures on $R_k$ with total variation distance.  We work with the sigma field $\F_H$ on $\Omega_H$ generated by all sets of the form $T_k^{-1}(\phi)$, where $\phi \in R_k$ and $k$ is a positive integer.  If $\mu$ is a measure on $\Omega_H$, it induces a measure $\mu \circ T_k^{-1}$ on $R_k$.  We claim that $\PP$ has a unique stationary distribution on the measure space $(\Omega_H, \F_H)$, a fact which will follow readily from the following lemma.

\begin{lemma}
	Let $W^1$ and $W^2$ be signed heaps starting from random initial configurations $W_0^1$ and $W_0^2$ with laws $\mu_0^1$ and $\mu_0^2$.  Then $\mu_n^1 = \PP^n \mu_0^1$ and $\mu_n^2 = \PP^n \mu_0^2$, the laws of $W_n^1$ and $W_n^2$ respectively, satisfy $\| \mu_n^1 \circ T_\ell^{-1} - \mu_n^2 \circ T_\ell^{-1} \| \rightarrow 0$ as $n \rightarrow \infty$ uniformly in $\mu_0^1$ and $\mu_0^2$.  Here $\| \cdot \|$ denotes total variation distance. \label{slemma}
\end{lemma}

\begin{proof}
	Construct a coupling between $W^1$ and $W^2$ as in Lemma \ref{basiclemma}.  Now suppose that at some point pieces fall successively in columns $0,1,2, \dots, m-1$ and that the pieces falling in column 0 did not annihilate.  Once this sequence is observed, the roofs of $W^1$ and $W^2$ must agree.  Call $\Phi_k$ the event that this sequence appears $k$ consecutive times.  In order for $T_\ell(W^1_n)$ and $T_\ell(W^2_n)$ to disagree after the event $\Phi_{k+\ell}$ is observed, at least $mk$ of the pieces comprising the $\Phi_{k+\ell}$ sequence must be annihilated.  By conditioning sequentially on these pieces and using Lemma \ref{blemma}, this probability is bounded above by $\left(\frac{1}{5}\right)^{mk}$.  Now, given $\epsilon>0$ we may choose $N,K$ large enough so that $\left(\frac{1}{5}\right)^{Km}<\frac{\epsilon}{2}$ and the probability that $\Phi_{K+\ell}$ is not observed in the first $N$ steps is less than $\frac{\epsilon}{2}$.  Since $\left\|\mu_n^1\circ T_\ell^{-1} - \mu_n^2\circ T_\ell^{-1} \right\|$ is bounded by the probability that $T_\ell(W^1_n)$ and $T_\ell(W^2_n)$ differ it follows that for all $n\geq N$, $\left\|\mu_n^1\circ T_\ell^{-1} - \mu_n^2\circ T_\ell^{-1} \right\|<\epsilon$.  
\end{proof}

Taking $\mu_0^1 = \mu_0$ and $\mu_0^2 = \mu_j$ the lemma implies that as $n \rightarrow \infty$, $\| \mu_{n}\circ T_\ell^{-1} - \mu_{n+j} \circ T_\ell^{-1} \| \rightarrow 0$ uniformly in $j$.  Thus $\mu_n \circ T_\ell^{-1}$ converges in $\MM(R_\ell)$ to a measure $\nu_\ell$, defined by $\nu_\ell(A) = \lim_{n \rightarrow \infty} \mu_n \circ T_\ell^{-1}(A)$ for all $A \in R_\ell$, which satisfies $\nu_\ell \circ T_\ell \circ T_k^{-1} = \nu_k$ for $k<\ell$.  It follows from Kolmogorov's extension theorem that there is a unique measure $\nu$ on $\Omega_H$ satisfying $\nu \circ T_k^{-1} = \nu_k$ for all $k$.  This measure is the unique stationary measure for $\PP$.  To see that $\PP \nu  = \nu$, observe that $\PP$ induces a continuous function $\tilde{\PP}:\MM(R_k) \rightarrow \MM(R_{k-1})$.  Since $\tilde{\PP}(\mu_n\circ T_k^{-1}) = \mu_{n+1} \circ T_{k-1}^{-1}$ and $\|\mu_n\circ T_k^{-1} - \nu_k \|,\|\mu_n\circ T_{k-1}^{-1} - \nu_{k-1} \| \rightarrow 0$ as $n \rightarrow \infty$, continuity implies that $\tilde{\PP}(\nu_k) = \nu_{k-1}$.  Hence, $(\PP \nu)\circ T_{k-1}^{-1} = \nu_{k-1}$ and it follows that $\PP \nu = \nu$.  Thus, $\nu$ is a stationary measure and by Lemma \ref{slemma} it is unique.  

\section{Concluding Remarks and Further Problems}
\begin{itemize}
\item[1.]We have shown that for a signed random heap, $\frac{1}{m}\lim_{n \rightarrow \infty} \E(\#T(W_n)) \leq 0.32893$.  However, the exact value of the limit is unknown for any $m \geq 4$.  

\item[2.]The random heap process may be generalized as follows.  Consider dropping pieces with sticky corners uniformly over the $m$ columns, and when a piece lands directly on top of another piece, flip a coin and allow the two pieces to annihilate with probability $p$.  When $p=0$, this new process is equivalent to the unsigned random heap process, and if $p = \frac{1}{2}$, we recover the signed random heap process.  Taking $p=1$, we obtain a random walk on the group 
\begin{equation}
\left< g_0,g_1, \dots, g_{m-1}:g_i g_j =g_j g_i \; \forall i \neq j \pm 1 \;\textrm{mod} \;m, \;\textrm{and}\; g_i^2 = e \;\forall i \right>. \nonumber
\end{equation}
  Numerical simulations for small values of $m$ suggest that the expected size of the roof is monotonic in $p$, but proving this remains a challenge.
\end{itemize}
\begin{acknowledgment}[Acknowledgements]
Many thanks to Yuval Peres for suggesting this problem and for helpful discussions and comments.  I also thank Bob Hough for helpful discussions, and G\'abor Pete for useful comments and permission to include his nice proof in the introduction.
\end{acknowledgment}



\end{document}